\documentclass[a4paper,reqno]{amsart}

\usepackage{amsmath,amssymb,amsthm}
\usepackage{fullpage}
\usepackage[colorlinks,linkcolor=blue,urlcolor=blue,citecolor=blue,destlabel=true]{hyperref}
\usepackage{graphicx}
\usepackage{stackrel}
\usepackage{tikz} 
\usetikzlibrary{matrix,arrows}
\usepackage{color}
\usepackage{textcomp}
\usepackage{enumerate}
\usepackage{mathtools}
\usepackage{subcaption}
\usepackage{algorithmic}
\usepackage[ruled,vlined,linesnumbered]{algorithm2e}
\usepackage{wrapfig}

\newtheoremstyle{myremark} 
    {7pt}                    
    {7pt}                    
    {}  	                 
    {}                           
    {\bf}       	         
    {.}                          
    {.5em}                       
    {}  
    
\theoremstyle{plain}
\newtheorem{lemma}{Lemma}[section]

\newtheorem{proposition}[lemma]{Proposition}
\newtheorem{theorem}[lemma]{Theorem}
\theoremstyle{definition}
\newtheorem{conjecture}[lemma]{Conjecture}
\newtheorem{definition}[lemma]{Definition}
\newtheorem{example}[lemma]{Example}
\newtheorem{question}[lemma]{Question}
\theoremstyle{myremark}

\newtheorem*{theorem*}{Theorem}

\newcommand{\bA}{\ensuremath{\mathbf{A}}}
\newcommand{\bB}{\ensuremath{\mathbf{B}}}
\newcommand{\bD}{\ensuremath{\mathbf{D}}}

\newcommand{\bGamma}{\ensuremath{\mathbf{\Gamma}}}

\newcommand{\bLambda}{\ensuremath{\mathbf{\Lambda}}}
\newcommand{\bX}{\ensuremath{\mathbf{X}}}
\newcommand{\C}{\ensuremath{\mathbb{C}}}

\newcommand{\N}{\ensuremath{\mathbb{N}}}

\newcommand{\R}{\ensuremath{\mathbb{R}}}

\newcommand{\cH}{\ensuremath{\mathcal{H}}}

\newcommand{\supp}{\ensuremath{\mathrm{supp}}}
\newcommand{\strain}{\ensuremath{\mathrm{Strain}}}

\newcommand{\tr}{\ensuremath{\mathrm{Tr}}}
\newcommand*\diff{\mathop{}\!\mathrm{d}}

\newcommand{\norm}[1]{\left\lVert#1\right\rVert}
\begin{document}
\author{Henry Adams}
\email{adams@math.colostate.edu}
\author{Mark Blumstein}
\email{blumstei@math.colostate.edu}
\author{Lara Kassab}
\email{kassab@math.colostate.edu}
\title{Multidimensional scaling on metric measure spaces}

\begin{abstract}
Multidimensional scaling (MDS) is a popular technique for mapping a finite metric space into a low-dimensional Euclidean space in a way that best preserves pairwise distances.
We overview the theory of classical MDS, along with its optimality properties and goodness of fit.
Further, we present a notion of MDS on infinite metric measure spaces that generalizes these optimality properties.
As a consequence we can study the MDS embeddings of the geodesic circle $S^1$ into $\mathbb{R}^m$ for all $m$, and ask questions about the MDS embeddings of the geodesic $n$-spheres $S^n$ into $\mathbb{R}^m$.
Finally, we address questions on convergence of MDS.
For instance, if a sequence of metric measure spaces converges to a fixed metric measure space $X$, then in what sense do the MDS embeddings of these spaces converge to the MDS embedding of $X$?
\end{abstract}

\maketitle

\setcounter{tocdepth}{1}

\section{Introduction}
\label{sec:intro}
Given $n$ objects and a notion of \emph{dissimilarity} between them, the classical multidimensional scaling (MDS) algorithm extracts a configuration of $n$ points in Euclidean space  whose pairwise distances ``best" approximate the given dissimilarities.
A typical source of dissimilarity data is the distance between high-dimensional objects, in which case MDS serves as a non-linear dimensionality reduction and visualization technique.
As such, the MDS algorithm is a popular technique for pattern recognition problems.
In this paper, we survey the classical algorithm, and describe an extension to (possibly infinite) metric measure spaces.
 
The coordinates extracted from an MDS embedding satisfy a least squares optimization problem.
While there are several popular choices of MDS loss function (metric or non-metric), we primarily focus on the classical algorithm which minimizes a form of loss function known as \emph{strain}.
The classical algorithm is algebraic and not iterative, simple to implement, and guaranteed to discover a configuration which optimizes the strain function.
Furthermore, if the input dissimilarities can be realized as distances in a Euclidean space, then classical MDS is guaranteed to recover such a configuration (unique up to translation and orthogonal transformation).
However, not all dissimilarity data admits a Euclidean realization.
In this case MDS produces a mapping into Euclidean space that distorts the inter-point pairwise distances as little as possible.
We make these ideas precise in Section~\ref{sec: MDS theory}.

The classical story is told using finite samples of points, finite dissimilarity matrices, and finite embedding coordinates.
Our goal is to extend to an infinite setting, where our input dissimilarity data is replaced by a \emph{metric measure space}: a metric space (with possibly infinitely many points) equipped with some probability measure.
This allows us to consider spaces whose points are weighted unequally, along with notions of convergence as more and more points are sampled from an infinite shape.

In more detail, a \emph{metric measure space} is a triple $(X,d,\mu)$ where $(X,d)$ is a compact metric space, and $\mu$ is a Borel probability measure on $X$.
In Section~\ref{sec: iMDS} we generalize the the classical MDS algorithm to metric measure spaces, and we show that this generalization minimizes the infinite analogue of strain.
As a motivating example, we consider the MDS embedding of the circle with the (non-Euclidean) geodesic metric, and equipped with the uniform measure.
By using the properties of circulant matrices, we identify the MDS embeddings of evenly-spaced points from the geodesic circle into $\R^m$, for all $m$.
As the number of points tends to infinity, these embeddings lie along the curve
\[\sqrt{2}\left(\cos\theta, \sin\theta, \tfrac{1}{3}\cos3\theta, \tfrac{1}{3}\sin3\theta,
\tfrac{1}{5}\cos5\theta, \tfrac{1}{5}\sin5\theta,\ldots\right)\in\R^m.\]
As this example illustrates, it is useful to consider the situation where a sequence of metric measure spaces $X_n$ converges to a fixed metric measure space $X$ as $n \to \infty$.
We survey various notions of convergence in Section~\ref{sec: convergence}.

Convergence is well-understood when each metric space has the same finite number of points, for example by Sibson's perturbation analysis~\cite{sibson1979studies}.
However, we are also interested in convergence when the number of points varies and is possibly infinite.
We survey results of~\cite{bengio2004learning,koltchinskii2000random} on the convergence of MDS when $n$ points $\{x_1,\ldots,x_n\}$ are sampled from a metric space according to a probability measure $\mu$, in the limit as $n\to\infty$.
The law of large numbers describes how the finite measures $\frac{1}{n}\sum_{i=1}^{n} \delta_{x_i}$ converge to $\mu$ as $n\to\infty$.
In~\cite{kassab2019multidimensional}, we reprove these results when instead we are given an \emph{arbitrary} sequence of probability measures $\mu_n\to\mu$.
The measures $\mu_n$ may now be unequally weighted, or have infinite support, for example.

\subsection*{Organization}
We present an overview on the theory of classical MDS in Section~\ref{sec: MDS theory}.
In Section~\ref{sec: Operator Theory}, we present necessary background information on operator theory and infinite-dimensional linear algebra.
We define a notion of MDS for infinite metric measure spaces in Section~\ref{sec: iMDS}.
In Section~\ref{sec: MDS circle}, we identify the MDS embeddings of the geodesic circle into $\R^m$, for all $m$, as a motivating example.
Lastly, in Section~\ref{sec: convergence}, we survey different notions of convergence of MDS.

\subsection*{Related Work}
\label{sec:related work}

The reader is referred to the introduction of~\cite{trosset1998new} and to~\cite{de198213, groenen2014past} for some aspects of the history of MDS.
There are a variety of papers that study some notion of robustness or convergence of MDS, including~\cite{bengio2004learning,sibson1978studies,sibson1979studies,sibson1981studies}.
Furthermore,~\cite[Section 3.3]{pekalska2001generalized} considers embedding new points in psuedo-Euclidean spaces,~\cite[Section 3]{diaconis2008horseshoes} considers infinite MDS in the case where the underlying space is an interval (equipped with some metric), and~\cite[Section 6.3]{buja2008data} discusses MDS on large numbers of objects.

\section{Classical Scaling}\label{sec: MDS theory}
Multidimensional scaling (MDS) is a set of statistical techniques concerned with the problem of using information about the dissimilarities between $n$ objects in order to construct a configuration of $n$ points in Euclidean space.
The input dissimilarities between the objects need not be based on Euclidean distances.

\begin{definition}
An $(n \times n)$ matrix $\mathbf D$ is called a \emph{dissimilarity matrix} if it is
symmetric and
\[d_{rr} = 0, \quad\text{with}\quad d_{rs} \geq 0 \quad\text{for all}\quad r \neq s.\]
\end{definition}
The first property above is called refectivity (the dissimilarity between an object and itself is zero), and the second property is called nonnegativity.
Symmetry requires that the dissimilarity from object $r$ to $s$ is the same as that from $s$ to $r$.
Note that there is no need to satisfy the triangle inequality.
A dissimilarity matrix $\mathbf D$ is called \emph{Euclidean} if there exists a
configuration of points in some Euclidean space whose interpoint distances are given by $\mathbf D$.

The goal of MDS is to map the objects $x_1,  \ldots, x_n$ to a configuration (or embedding) of points $f(x_1), \ldots, f(x_n)$ in $\R^m$ so that the given dissimilarities $d(x_i,x_j)$ are well-approximated by the Euclidean distances $\|f(x_i) - f(x_j)\|_2$.
The different notions of approximation give rise to the different types of MDS.

If the dissimilarity matrix can be realized exactly as the distance matrix of some set of points in $\R^m$ (i.e.\ if the dissimilarity matrix is \emph{Euclidean}), then MDS will find such a realization.
Furthermore, MDS can be used to identify the minimum Euclidean dimension $m$ admitting such an isometric embedding.
However, some dissimilarity matrices or metric spaces are inherently non-Euclidean (cannot be embedded into $\R^m$ for any $m$).
When a dissimilarity matrix is not Euclidean, then MDS produces a mapping into $\R^m$  that distorts the interpoint pairwise distances as little as possible.
Though we introduce MDS below, the reader is also referred to~\cite{bibby1979multivariate, cox2000multidimensional, groenen2014past} for more complete introductions.

Classical multidimensional scaling (cMDS) is also known as Principal Coordinates Analysis (PCoA), Torgerson Scaling, or Torgerson--Gower scaling.
The cMDS algorithm minimizes a loss function called \emph{strain}, and one of the main advantages of cMDS is that its algorithm is algebraic and not iterative.
Therefore, it is simple to implement, and it is guaranteed to discover the optimal configuration in $\R^m$.
In this section, we describe the cMDS algorithm, and then discuss some of its optimality properties and goodness of fit.

\begin{figure}[htb]
\centering
\includegraphics[scale=0.3]{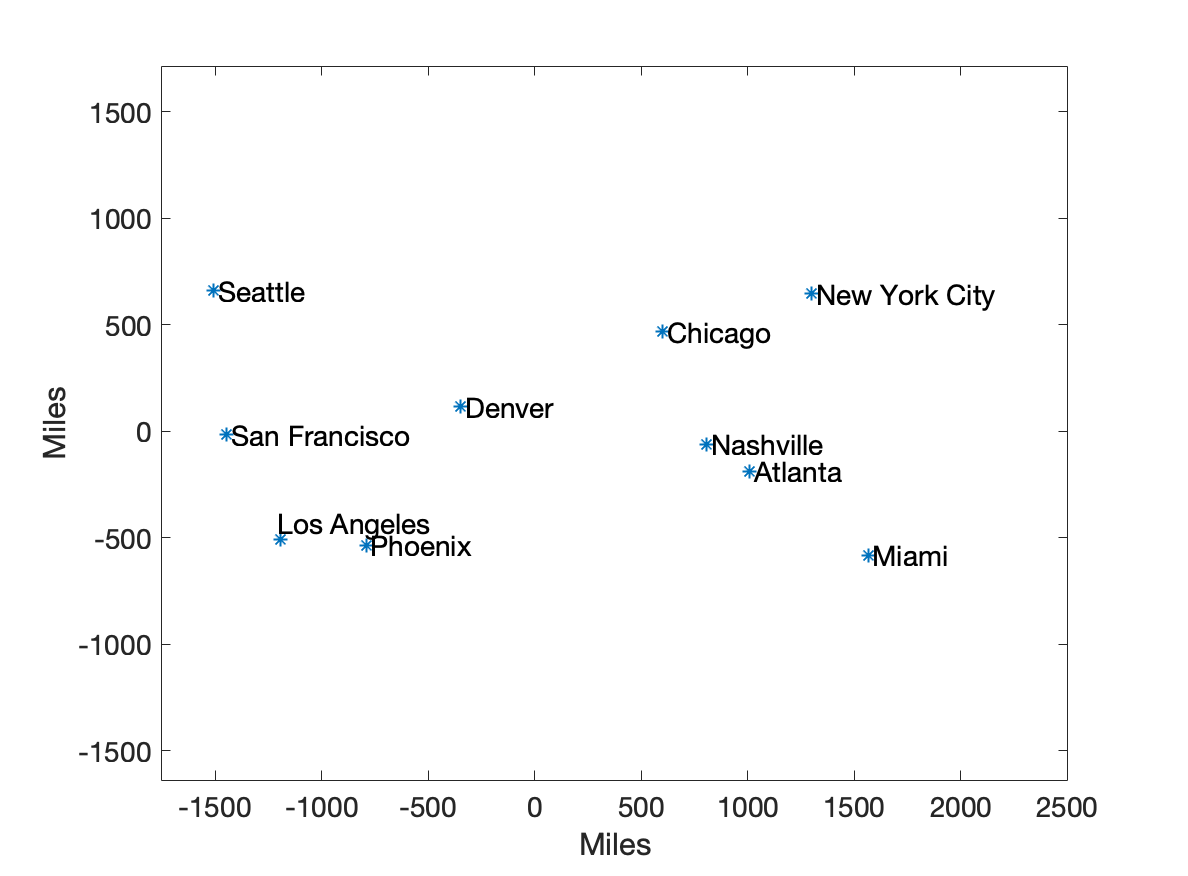}
\caption{Reconstruction of a map of 10 U.S.\ cities using MDS}
\label{fig: MDS cities}
\end{figure}

As an illustrative example, we consider ten U.S.\ cities equipped with the road distance between them, which is a non-Euclidean distance.
The classical MDS algorithm produces a two dimensional configuration of points (see Figure~\ref{fig: MDS cities}), where the points represent the different cities.
The Euclidean pairwise distances (distances as the crow flies) between the cities in the MDS embedding are the Euclidean distances that best approximate the road distances between them.

Let $\bD = (d_{ij})$ be an $n\times n$ dissimilarity matrix.
Let $\bA = (a_{ij})$, where $a_{ij} = -\frac{1}{2}d^2_{ij}$.
Define the matrix $\bB$ to be the double mean-centering of $\bA$, with entries given by
\begin{equation}\label{eq:B-definition}
b_{rs}=a_{rs}-\frac{1}{n}\sum\limits_{s=1}^n a_{rs}-\frac{1}{n}\sum\limits_{r=1}^n a_{rs}+\frac{1}{n^2}\sum\limits_{r,s=1}^n a_{rs}.
\end{equation}

 Since $\bD$ is a symmetric matrix, it follows that $\bA$ and $\bB$ are each symmetric, and therefore $\bB$ has $n$ real eigenvalues.

Assume for convenience that there are at least $m$ positive eigenvalues for matrix $\bB$, where $m\le n$.
By the spectral theorem of symmetric matrices, let $\bB = \bGamma \bLambda \bGamma^\top$, with $\bGamma$ containing unit-length eigenvectors of $\bB$ as its columns, and with the diagonal matrix $\bLambda$ containing the eigenvalues of $\bB$ in decreasing order along its diagonal.
Let $\bLambda_m$ be the $m\times m$ diagonal matrix of the largest $m$ eigenvalues sorted in descending order, and let $\bGamma_m$ be the $n\times m$ matrix of the corresponding $m$ eigenvectors in $\bGamma$.
The coordinates of the MDS embedding into $\R^m$ are then given by the $n\times m$ matrix $\bX=\bGamma_m\bLambda_m^{1/2}$.
More precisely, the MDS embedding consists of the $n$ points in $\R^m$ given by the $n$ rows of $\bX$.
The procedure for classical MDS can be summarized in the following algorithm.

\begin{algorithm}[H]
\SetAlgoLined
  \SetKwInOut{Input}{input}
    \Input{Dissimilarity matrix $\bD$}
     \SetKwInOut{Compute}{compute}
     \Compute{Compute the matrix $\bA = (a_{ij})$, where $a_{ij} = -\frac{1}{2}d^2_{ij}$}
     \SetKwInOut{Perform}{perform}
     \Perform{Perform double-centering to $\bA$: define $\bB$ by~\eqref{eq:B-definition}}
     \Compute {Compute the eigendecomposition of $\bB = \bGamma \bLambda \bGamma^\top$}
      \SetKwInOut{Select}{select}
      \Select {Select the $m$ largest nonnegative eigenvalues of $B$ to obtain $ \bLambda_m$}
     \SetKwInOut{Output}{output}
     \Output{Coordinate matrix given by $\bX=\bGamma_m\bLambda_m^{1/2}$}
    
 \caption{Classical MDS}
 \label{alg: cMDS}
\end{algorithm}

We give a small example.

\begin{example}
\label{ex: non-Euclidean}
We implement Algorithm~\ref{alg: cMDS} on the following $4\times 4$ dissimilarity matrix $\bD$.
\[
\bD=\begin{pmatrix}
0&2&2&1\\
2&0&2&1\\
0&2&2&1\\
1&1&1&0
\end{pmatrix}
\quad
\bA=-\frac{1}{2}\begin{pmatrix}
0&4&4&1\\
4&0&4&1\\
0&4&4&1\\
1&1&1&0
\end{pmatrix}
\quad
\bB=\frac{1}{16}\begin{pmatrix}
\phantom{-}13&-15&\phantom{-}5&-3\\
-15&\phantom{-}21&-7&\phantom{-}1\\
\phantom{-}5&-7&-3&\phantom{-}5\\
-3&\phantom{-}1&\phantom{-}5&-3
\end{pmatrix}
\]
The eigenvalues of $\bB$ are $2.159\ldots$, $0.192\ldots$, $0$, and $-0.602\ldots$, and the MDS embedding of $\bD$ in $\R^2$ is drawn in Figure~\ref{fig: non-Euclidean}.

\begin{figure}[htb]
\centering
\def\svgwidth{0.2\linewidth}
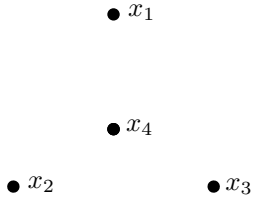
\caption{MDS embedding of the $4\times 4$ distance matrix in Example~\ref{ex: non-Euclidean}.}
\label{fig: non-Euclidean}
\end{figure}

This dissimilarity matrix is not Euclidean.
Indeed, label the points $x_1$, $x_2$, $x_3$, $x_4$ in order of their row/column in $\bD$.
In any isometric embedding in $\R^n$, the points $x_1$, $x_2$, $x_3$ would be mapped to an equilateral triangle.
The point $x_4$ would need to get mapped to the midpoint of each edge of this triangle, which is impossible in Euclidean space.
\end{example}

The following fundamental criterion determines algebraically whether a dissimilarity matrix $\bD$ is Euclidean or not.\\

\begin{theorem}~\cite[Theorem~14.2.1]{bibby1979multivariate}\label{thm: MDS p.s.d}
Let $\bD$ be a dissimilarity matrix, and define $\bB$ by equation~\eqref{eq:B-definition}.
Then $\bD$ is Euclidean if and only if $\bB$ is a positive semi-definite matrix.
\end{theorem}

Moreover, if $\bB$ is positive semi-definite of rank $m$, then a perfect realization of the dissimilarities can be found by a collection of points in $m$-dimensional Euclidean space.

Let $\bD$ be a dissimilarity matrix, and define $\bB$ via~\eqref{eq:B-definition}.
A measure of the goodness of fit of MDS, even in the case when $\bD$ is not Euclidean, can be obtained as follows.
If $\hat{\bX}$ is a fitted configuration in $\R^m$ with centered inner-product matrix $\hat{\bB}$, then a measure of the discrepancy between $\bB$ and $\hat{\bB}$ is the following strain function~\cite{mardia1978some},
\begin{equation}\label{eq:optimization}
\mathrm{tr}((\bB-\hat{\bB})^2)=\sum\limits_{i,j=1}^n(b_{i,j}-\hat{b}_{i,j})^2.
\end{equation}

\begin{theorem}~\cite[Theorem~14.4.2]{bibby1979multivariate}\label{thm: strain-minimization-cMDS}
Let $\bD$ be a dissimilarity matrix.
Then for fixed $m$, the strain function in \eqref{eq:optimization} is minimized over all configurations $\hat {\bX}$ in $m$
dimensions when $\hat {\bX}$ is the classical solution to the MDS problem.
\end{theorem}

The reader is referred to~\cite[Section~2.4]{cox2000multidimensional} for a summary of a related optimization procedure with a different normalization, due to Sammon~\cite{sammon1969nonlinear}.

\section{Preliminaries}
\label{sec: Operator Theory}

\begin{figure}[htb]
\centering
\begin{minipage}{0.49\linewidth}
\centering
\includegraphics[width=1\textwidth]{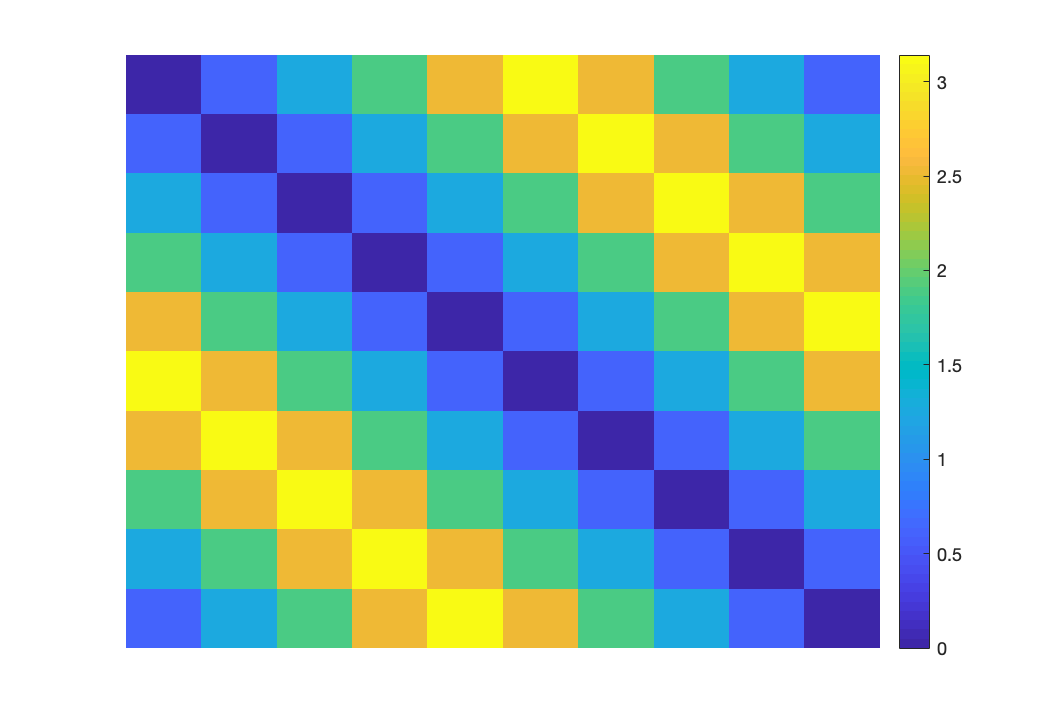}
 \subcaption{Distance matrix of ten points.}
\label{fig:sub:distancematrix}
\end{minipage}
\begin{minipage}{0.49\linewidth}
\centering
\includegraphics[width=1\textwidth]{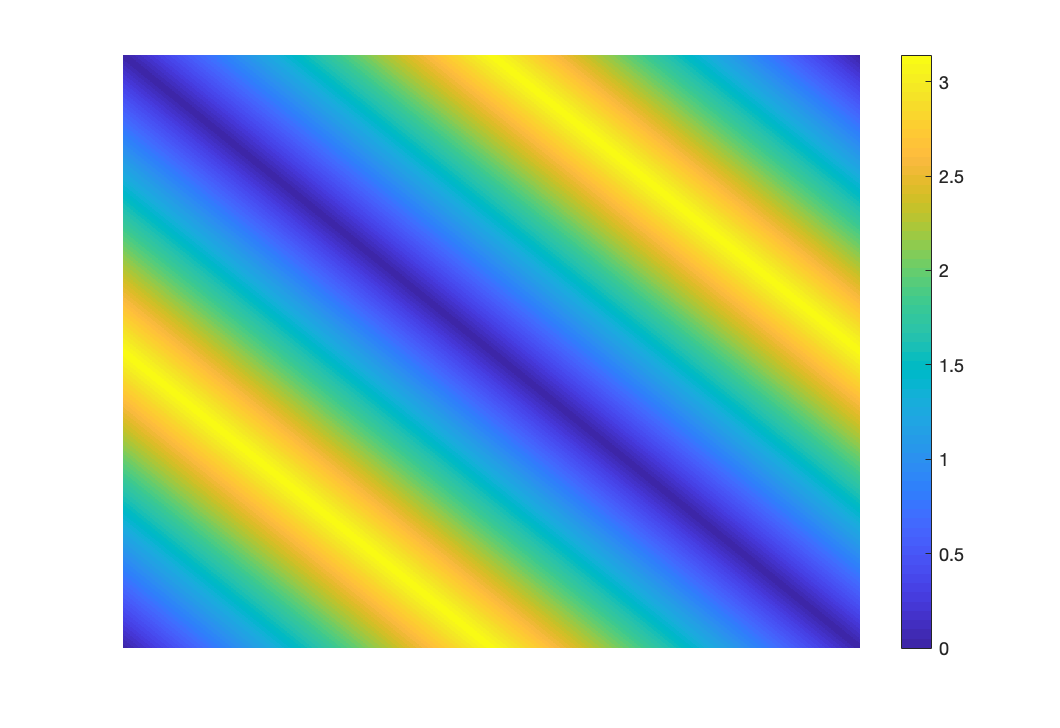}
 \subcaption{Distance kernel of infinitely many points.}
\label{fig:sub:distancekernel}
\end{minipage}
      \caption{Distance matrix or kernel corresponding to (A) ten equally-spaced points, and to (B) all of the points, on the circle equipped with the geodesic metric.
      }

    \label{fig: matrix-to-kernel}
\end{figure}

We are interested in studying the MDS embeddings of spaces with possibly infinitely many points, and distance matrices aren't enough to store infinitely many pairwise distances.
Instead, we use \emph{kernels}, which roughly speaking are distance functions that compute the pairwise distance between any two points in the space.
For example, the kernel corresponding to the geodesic distance on a circle is illustrated in Figure~\ref{fig:sub:distancekernel}.

This section introduces the reader to concepts in infinite-dimensional linear algebra and operator theory used throughout the paper.

\subsection*{Kernels and Operators}
Let $X$ be a metric space equipped with a measure $\mu$.
We denote by $L^2(X, \mu)$  the set of square-integrable real-valued $L^2$-functions with respect to the measure $\mu$.
We note that $L^2(X, \mu)$ is furthermore a Hilbert space, after equipping it with the inner product given by $\langle f,g \rangle = \int_X fg\ d\mu$.

A \emph{real-valued $L^2$-kernel} $K \colon X \times X \to \R$ is a continuous measurable square-integrable function.
The kernels that we consider in this paper are symmetric, meaning that $K(x,s) = K(s,x)$ for all $x,s\in X$.
A symmetric kernel is \emph{positive semi-definite} if
\[{\sum\limits _{i=1}^{n}\sum\limits _{j=1}^{n}c_{i}c_{j}K(x_{i},x_{j})\geq 0}\]
holds for any $m\in \N$, any $x_1,\dots ,x_m\in X$, and any $c_1,\dots,c_m\in \R$.
At least in the case when $X$ is a compact subspace of $\R^m$ (and probably more generally), a symmetric kernel is positive semi-definite if 
\[ \int \limits _{X} \int\limits _{X} K(x, s)f(x)f(s)\mu(\diff x)\mu (\diff s) \geq 0 \]
for any $f \in L^2(X, \mu)$.

\begin{definition}[Hilbert--Schmidt Integral Operator]\label{def: Kernel Operator}
Let $(X, \Omega, \mu)$ be a $\sigma$-finite measure
space, and let $K$ be an $L^2$-kernel on $X \times X$.
Then the integral operator
\[ [T_K\phi](x) = \int \limits_X K(x,s) \phi(s) \mu(\diff s),\]
which defines a linear mapping from the space $L^2 (X, \mu)$ into itself, is called a \emph{Hilbert--Schmidt integral operator}.
\end{definition}

Hilbert--Schmidt integral operators are both continuous (and hence bounded) and compact operators.

\begin{definition}
A Hilbert--Schmidt integral operator is a \emph{self-adjoint operator} if $K(x,y) = {K(y,x)}$  holds for almost all $(x,y) \in X \times X$ (with respect to $\mu \times \mu$).
\end{definition}

\begin{definition}\label{p.s.d.o}
A bounded self-adjoint operator $T$ on a Hilbert space $\cH$ is a \emph{positive semi-definite operator} if $\langle Tx, x \rangle \geq 0$
for any $x \in \cH$.
\end{definition}

It follows that the eigenvalues of a positive semi-definite operator $A$, when they exist, are real.

\subsection*{The Spectral Theorem}

Classical MDS relies on the fact that symmmetric matrices are orthogonally diagonalizable with real eigenvalues.
Furthermore, positive semi-definite matrices (having nonnegative eigenvalues) may be represented as matrices of Euclidean inner products. The following two theorems give analogues of these results for kernels instead of matrices.

\begin{theorem}[Spectral theorem on compact self-adjoint operators]
\label{thm:spectral}
Let $\cH$ be a 
Hilbert space, and suppose $T \colon \cH \to \cH$ is a bounded compact self-adjoint operator.
Then $T$ has at most a countable number of nonzero eigenvalues $\lambda_n  \in \R$, with a corresponding orthonormal set $\{ e_n \}$ of eigenvectors, such that
\[ T(\cdot) = \sum\limits_{n=1}^\infty \lambda_n \langle e_n, \cdot \rangle\ e_n.\]
Furthermore, the multiplicity of each nonzero eigenvalue is finite, zero is the only possible accumulation point of $\{ \lambda_n \}$, and if the set of nonzero eigenvalues is infinite then zero is an accumulation point.
\end{theorem}

A fundamental theorem that characterizes positive semi-definite kernels is the Generalized Mercer's Theorem.\\

\begin{theorem}[Generalized Mercer's Theorem]~\cite[Lemma 1]{kuhn1987eigenvalues}\label{Thm: Mercer's Theorem}
Let $X$ be a compact topological Hausdorff space equipped with a finite Borel measure $\mu$, and let $K \colon X \times  X \to \C$ be a continuous positive semi-definite kernel.
Then, there exists a scalar sequence $\{ \lambda_n \} \in \ell_1$ with $\lambda_1 \geq \lambda_2 \geq \cdots \geq 0$, and an orthonormal system $\{ \phi_n\}$ of continuous square-integrable functions with respesct to $\mu$,
such that the expansion 
\begin{equation*}\label{eq: Mercer}
K(x, s) = \sum\limits _{n=1} ^\infty \lambda_n\phi_n(x)\overline {\phi_n(s)}, \quad x, s \in \supp(\mu)
\end{equation*} converges uniformly, where $\supp$ denotes the support of a measure $\mu$.
\end{theorem} 

Therefore, given $X$ and $K$ as in Theorem~\ref{Thm: Mercer's Theorem}, the associated Hilbert--Schmidt integral operator
\begin{equation*}\label{eq: Kernel Operator}
[T_K\phi](x) = \int \limits_X K(x,s) \phi(s) \mu(\diff s)\end{equation*}
is also positive semi-definite.
Moreover, the eigenvalues of $T_K$ can be arranged in non-increasing order $\lambda_1 \geq \lambda_2 \geq \ldots \geq 0$, indexed according to their algebraic multiplicities, and the orthonormal system $\{ \phi_n \}$ gives the corresponding eigenfunctions of $T_K$.

\section{MDS of Infinite Metric Measure Spaces}
\label{sec: iMDS}

Classical multidimensional scaling (cMDS) can be described either as a strain-minimization problem, or as a linear algebra algorithm involving eigenvalues and eigenvectors.
Indeed, one of the main theoretical results for cMDS is that the linear algebra algorithm solves the corresponding strain-minimization problem (see Theorem~\ref{thm: strain-minimization-cMDS}).
In this section, we describe how to generalize both of these formulations to (possibly infinite) metric measure spaces.

This will allow us to discuss the MDS embedding of the circle, for example, without needing to restrict attention to finite subsets thereof.

\begin{definition} A \emph{metric measure space} is a triple $(X,d,\mu)$ where
\begin{itemize}
\item $(X,d)$ is a compact metric space, and
\item $\mu$ is a Borel probability measure on $X$, i.e.\ $\mu (X) = 1$.
\end{itemize}
\end{definition}

Given a metric space $(X, d)$, by a measure on $X$ we mean a measure on the Borel $\sigma$-algebra of $X$.
When it is clear from the context, the triple $(X, d, \mu)$ will be denoted by only
$X$.
The reader is referred to~\cite{memoli2011gromov, memoli2014gromov} for details on metric measure spaces, and for interpretations of these concepts in the context of object matching.

Let $(X, d, \mu)$ be a metric measure space, with $d$ a $L^2$-function on $X \times X$.
We say that $X$ is \emph{Euclidean} if it can be isometrically embedded into $(\ell^2,\| \cdot \|_2)$.
$X$ is furthermore Euclidean in the finite-dimensional sense if there is an isometric embedding $X\to \R^m$.

\subsection*{MDS on Infinite Metric Measure Spaces}\label{ss:infinite-mds}
Let $(X, d, \mu)$ be a metric measure space, where $d$ is an $L^2$-function on $X \times X$.

We propose the following MDS method on infinite metric measure spaces:

\begin{enumerate}[(i)]
\item From the metric $d$, construct the kernel $K_A \colon X \times X \to \R$ defined as 
\begin{equation}\label{eq : kernel A} K_A(x,s)= -\tfrac{1}{2}d^2(x,s).
\end{equation}
\item Obtain the kernel $K_B \colon X \times X \to \R$ via
\begin{equation}\label{eq : kernel B}
K_B(x,s)= K_A(x,s) - \int \limits_X  K_A(w,s) \mu(\diff w) - \int \limits_X  K_A(x,z) \mu(\diff z) + \int\limits_{X} \int\limits_X K_A(w,z) \mu(\diff w) \mu(\diff z).
\end{equation}
Assume $K_B \in L^2(X \times X)$.
Define $T_{K_B}\colon L^2(X) \to L^2(X)$ as
\begin{equation}\label{eq: operator TB}
[T_{K_B}\phi](x) = \int \limits_X K_B(x,s) \phi(s) \mu(\diff s).
\end{equation}

\item Let $\lambda_1\geq \lambda_2 \geq \dots$ denote the eigenvalues of $T_{K_B}$, with corresponding eigenfunctions $\phi_1, \phi_2, \ldots \in L^2(X)$ forming an orthonormal system in $L^2(X)$.

\item Define $K_{\hat B}(x,s) = \sum\limits_{i=1} ^ \infty \hat \lambda_i \phi _i (x) \phi_i (s)$, where $\hat \lambda_i=\lambda_i$ if $\lambda_i\geq 0$, and otherwise $\hat \lambda_i=0$.
Let $T_{K_{\hat B}}\colon L^2(X) \to L^2(X)$ be the Hilbert--Schmidt integral operator associated to the kernel $K_{\hat B}$.
The eigenfunctions $\phi _i$ for $T_{K_B}$ (with eigenvalues $\lambda_i$) are also the eigenfunctions for $T_{K_{\hat B}}$ (with eigenvalues $\hat \lambda_i$).
By Mercer's Theorem (Theorem~\ref{Thm: Mercer's Theorem}), $K_{\hat B}$ converges uniformly.

\item Define the MDS embedding of $X$ into $\ell^2$ via the map $f\colon X\to \ell^2$ given by 
\begin{equation}\label{eq: MDS embedding}
    f(x)=\left(\sqrt{\hat \lambda_1} \phi_1(x),\sqrt{\hat \lambda_2} \phi_2(x),\sqrt{\hat \lambda_3} \phi_3(x),\ldots\right)\end{equation}
Similarly, define the MDS embedding of $X$ into $\R^m$ via the map $f_m\colon X\to \R^m$ given by
\[f_m(x)=\left(\sqrt{\hat \lambda_1} \phi_1(x), \sqrt{\hat \lambda_2} \phi_2(x), \ldots, \sqrt{\hat \lambda_m} \phi_m(x)\right)\]
\end{enumerate} 

The procedure for infinite classical MDS can be summarized in the following algorithm.

\begin{algorithm}[H]
\SetAlgoLined
  \SetKwInOut{Input}{input}
    \Input{Metric measure space $(X, d, \mu)$}
     \SetKwInOut{Compute}{compute}
     \Compute{Compute the kernel $K_A$ by~\eqref{eq : kernel A}}
     \SetKwInOut{Perform}{perform}
     \Perform{Perform double-centering to $K_A$: define the kernel $K_B$ by~\eqref{eq : kernel B}}
     \Compute {Compute the eigenvalues and eigenfunctions of the operator $T_{K_B}$ defined by~\eqref{eq: operator TB}}
     \SetKwInOut{Select}{select}
      \Select {Select the nonnegative eigenvalues of $T_{K_B}$ to obtain a new operator $T_{K_{\hat B}}$}
     \SetKwInOut{Output}{output}
     \Output{MDS embedding of $X$ into $\ell^2$ via the map $f\colon X\to \ell^2$ given in~\eqref{eq: MDS embedding}}
    
 \caption{Classical MDS on Infinite Spaces}
 \label{alg: infinite cMDS}
\end{algorithm}

\begin{table}[h]
    \caption{A comparison of various aspects of classical and infinite MDS.
This table is constructed analogously to that on~\cite{wiki:FPCA} Principal Component Analysis (PCA) and Functional Principal Component Analysis (FPCA).
    }
    \label{table:sample}
    \begin{center}
\renewcommand{\arraystretch}{2.7}
\resizebox{\textwidth}{!}{
\begin{tabular}{ |p{4cm}||p{5.5cm}|p{6.5cm}|  }
 \hline
 \textbf{Elements} &  \textbf{Classical MDS}  &  \textbf{Infinite MDS}\\
 \hline
Data  & $(X,d)$ with $|X|<\infty$    & $(X, d, \mu)$ \\
 \hline
Distance Representation &   $D_{i,j} = d(x_i, x_j),\quad D \in \mathcal M _{n \times n}$ & $K_D(x,s) = d(x,s) \in L^2_{\mu \otimes \mu }(X\times X)$  \\
 \hline
Linear Operator &
\begin{tabular}{@{}c@{}}$b_{rs}=a_{rs}-\frac{1}{n}\sum\limits_{s=1}^n a_{rs}$ \\ \hspace{10mm}$-\frac{1}{n}\sum\limits_{r=1}^n a_{rs}+\frac{1}{n^2}\sum\limits_{r,s=1}^n a_{rs}$\end{tabular}
& $[T_{K_B}\phi](x) = \int \limits_X K_B(x,s) \phi(s) \mu(\diff s)$ \\
  \hline
 Eigenvalues & $\lambda_1 \geq \lambda_2 \geq \ldots \geq \lambda_n $ & $\hat \lambda_1 \geq \hat \lambda_2 \geq \ldots $\\
  \hline
Eigenvectors &$v^{(1)}, v^{(2)}, \ldots, v^{(m)} \in \R^n$ & $\phi_1(x), \phi_2(x), \ldots \in L^2(X) $ \\
 \hline
Embedding in $\R^m$ or $\ell^2$& \footnotesize{$  f(x_i) = \left(\sqrt{ \lambda_1} v^{(i)}_1,\sqrt{ \lambda_2} v^{(i)}_2, \ldots, \sqrt{ \lambda_m} v^{(i)}_m \right)$}  &  \footnotesize{$f(x) = \left(\sqrt{\hat \lambda_1} \phi_1(x),\sqrt{\hat \lambda_2} \phi_2(x),\sqrt{\hat \lambda_3} \phi_3(x),\ldots\right)$}  \\
 \hline
 Strain Minimization& $\sum\limits_{i,j=1}^n(b_{i,j}-\hat{b}_{i,j})^2$  &$\int \int \left( K_B(x,s) - K_{\hat B}(x,s) \right)^2 \mu(\diff x) \mu(\diff s)$ \\
 \hline
\end{tabular}}
    \end{center}
\end{table}

\begin{proposition}\cite[Proposition 6.3.1.]{kassab2019multidimensional}
The MDS embedding map $f\colon X\to \ell^2$ 
is continuous.
\end{proposition}

The following theorem generalizes Theorem~\ref{thm: MDS p.s.d} to metric measure spaces.

\begin{theorem}\cite[Theorem 6.3.3.]{kassab2019multidimensional}\label{thm:inf-mds-euc}
A metric measure space $(X, d, \mu)$ is Euclidean if and only if $T_{K_B}$ is a positive semi-definite operator on $L^2(X, \mu).$
\end{theorem}

We show that MDS for metric measure spaces minimizes the loss function $\strain(f)$, defined as
\[\strain(f) =\|T_{K_B}-T_{K_{\hat{B}}}\|_{HS}^{2}=\tr ((T_{K_B}-T_{K_{\hat{B}}})^2) = \int \int \left( K_B(x,s) - K_{\hat B}(x,s) \right)^2 \mu(\diff x) \mu(\diff s).\]

This result generalizes~\cite[Theorem~14.4.2]{bibby1979multivariate}, or equivalently~\cite[Theorem~2]{trosset1997computing}, to the infinite case.

\begin{theorem}\cite[Theorem 6.4.3.]{kassab2019multidimensional}\label{Thm: infinite-mds-optimization}
Let $(X, d, \mu)$ be a metric measure space.
Then $\strain(f)$ is minimized over all maps $f\colon X\to \ell^2$ or $f\colon X\to \R^m$ when $f$ is the MDS embedding given in Section~\ref{ss:infinite-mds}.
\end{theorem}

\renewcommand{\arraystretch}{1}

\section{MDS of the Circle}
\label{sec: MDS circle}
Let $S^1$ be the unit circle equipped with arc-length distance and the uniform measure ($\frac{d\theta}{2\pi}$).
Using our definition of MDS as an integral operator, we show that MDS maps $S^1$ into an infinite dimensional sphere of radius $\frac{\pi}{2}$ sitting inside $\ell^2$.
The embedded circle occupies an infinite number of dimensions in $\ell^2$, and in fact, the infinite dimensional space is needed---the embedding is better (in the sense of strain minimization) than the MDS embedding into $\R^m$ for any finite $m$.

It is instructive to consider how MDS on finite samples of $S^1$ converges to the MDS integral operator on the entire circle.
We start with the easiest case: let $S^1_n$ be the sample of $n$ evenly-spaced points on $S^1$.
\begin{proposition}\label{prop:S1n}
 The classical MDS embedding of $S^1_n$ lies, up to a rigid motion of $\R^m$, on the curve $\gamma_m\colon S^1_n\to\R^m$ defined by
\[ \gamma_m(\theta) = (a_{1,n}\cos(\theta),a_{1,n}\sin(\theta),a_{3,n}\cos(3\theta),a_{3,n}\sin(3\theta),a_{5,n}\cos(5\theta),a_{5,n}\sin(5\theta),\ldots)\in\R^m, \]
where $\lim_{n\to\infty} a_{j,n}= \frac{\sqrt{2}}{j}$ (with $j$ odd).
\end{proposition}

\begin{figure}[htb]
\centering  \includegraphics[width=0.5\textwidth]{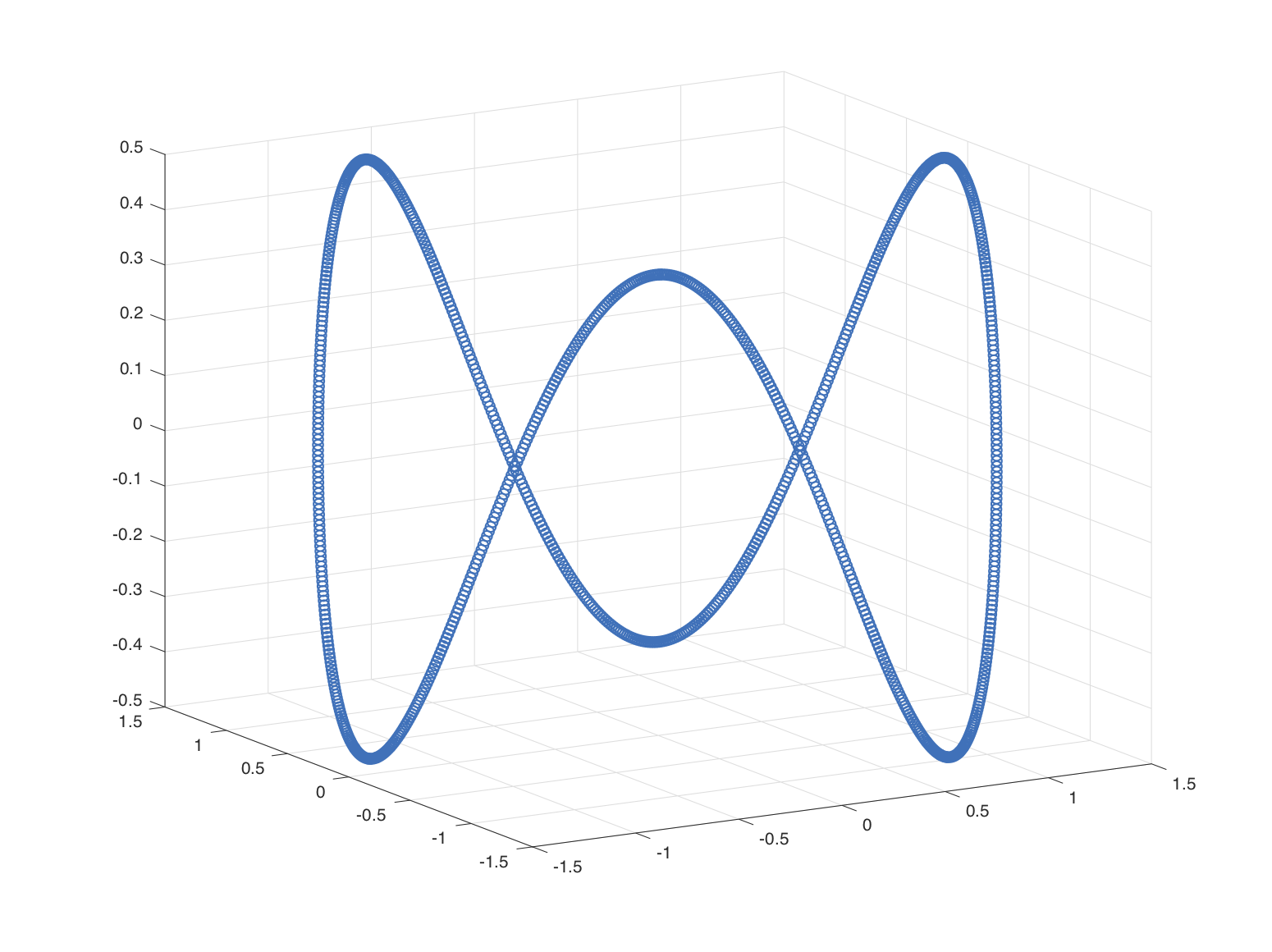}
\caption{MDS embedding of $S^1_{1000}$.}
\label{fig: MDS of Circle}
\end{figure}

Figure~\ref{fig: MDS of Circle} shows the MDS configuration in $\R ^3$ of 1000 equally-spaced points on $S^1$ obtained using the three largest positive eigenvalues.

We sketch the outline of this computation; full details are given in ~\cite{kassab2019multidimensional}.
Let $\bD$ be the arc-length distance matrix for $S^1_n$.
Following the steps of classical MDS, define $\bA = (a_{ij})$ with $a_{ij} = -\frac{1}{2}d^2_{ij}$, and let $\bB$ be the doubly mean-centered version of matrix $\bA$.
A matrix $\mathbf{M}$ is called \emph{circulant} if cyclically shifting all rows of $\mathbf{M}$ down by one has the same effect as cyclically shifting all columns of $\mathbf{M}$ left by one.
Both $\bD$ and the double mean centering matrix have this property, and therefore the MDS symmetric matrix $\bB$ is circulant.
In coordinates, it has the following form: \begin{equation*}\label{eq:symmetric-circulant}
\bB =
\begin{pmatrix}
b_{0} & b_{1} & b_{2} & \ldots & b_{3} & b_{2} & b_{1} \\
b_{1} & b_{0} & b_{1} & \ldots & b_{4} & b_{3} & b_{2} \\
b_{2} & b_{1} & b_{0} & \ldots & b_{5} & b_{4} & b_{3} \\
\vdots & \vdots & \vdots & & \vdots & \vdots & \vdots \\
b_{3} & b_{4} & b_{5} & \ldots & b_{0} & b_{1} & b_{2} \\
b_{2} & b_{3} & b_{4} & \ldots & b_{1} & b_{0} & b_{1} \\
b_{1} & b_{2} & b_{3} & \ldots & b_{2} & b_{1} & b_{0}
\end{pmatrix}.
\end{equation*}

For example, if $\bD$ is the distance matrix for $n=7$ equally-spaced points on the circle, then we compute
\[
\bD=\frac{2\pi}{7}\begin{pmatrix}
0&1&2&3&3&2&1\\
1&0&1&2&3&3&2\\
2&1&0&1&2&3&3\\
3&2&1&0&1&2&3\\
3&3&2&1&0&1&2\\
2&3&3&2&1&0&1\\
1&2&3&3&2&1&0
\end{pmatrix}
\quad
\text{and}
\quad
\bB=\frac{2\pi^2}{49}\begin{pmatrix}
\phantom{-}4&\phantom{-}3&\phantom{-}0&-5&-5&\phantom{-}0&\phantom{-}3\\
\phantom{-}3&\phantom{-}4&\phantom{-}3&\phantom{-}0&-5&-5&\phantom{-}0\\
\phantom{-}0&\phantom{-}3&\phantom{-}4&\phantom{-}3&\phantom{-}0&-5&-5\\
-5&\phantom{-}0&\phantom{-}3&\phantom{-}4&\phantom{-}3&\phantom{-}0&-5\\
-5&-5&\phantom{-}0&\phantom{-}3&\phantom{-}4&\phantom{-}3&\phantom{-}0\\
\phantom{-}0&-5&-5&\phantom{-}0&\phantom{-}3&\phantom{-}4&\phantom{-}3\\
\phantom{-}3&\phantom{-}0&-5&-5&\phantom{-}0&\phantom{-}3&\phantom{-}4
\end{pmatrix}.
\]

The complex eigenvectors of such a matrix are given by the \emph{discrete Fourier modes}, namely $ x_k(n) := (w_{n}^0, w_{n}^k, \ldots, w_{n}^{(n-1)k})^\top$ for $0\le k\le n-1$, where $ w_{n} = e^{2\pi i/n}$.
Since the first entry of each vector $x_k$ is one, the eigenvalue of $x_k$ can be computed simply by taking the dot product of the first row of $\bB$ with $x_k$.
Note that the vector of all ones has eigenvalue zero.

Since $\bB$ is symmetric, each complex eigenvector can be split into its real and imaginary part, which forms two real eigenvectors---this explains the sine and cosine representation of eigenvectors in the proposition.
It turns out that the odd Fourier modes have positive eigenvalues, and the even Fourier modes have negative eigenvalues.
Since MDS retains coordinates corresponding to positve eigenvalues, we are left with only the odd Fourier modes.

How does this finite MDS computation compare to the MDS integral operator on all of $S^1$? Let $S^1$ be the unit circle with arc-length distance and uniform measure.
If $\phi_k(x)= e^{ikx}$, then one may check (use integration by parts) that
\begin{align*}
-\frac{1}{2} \int_{y=x-\pi}^{y=x+\pi} (y-x)^2 e^{iky} \; \frac{dy}{2\pi} &=\frac{1}{k^2}(-1)^{k+1} e^{ikx}.
\end{align*} 
Despite not having performed the double mean centering step to the kernel function, this computation shows that the (complex) eigenfunctions of MDS on $S^1$ are $\phi_k(x)= e^{ikx}$ with $\lambda_k=\frac{1}{k^2}(-1)^{k+1}$, $k \neq 0$.
Indeed, the mean centering step associates the eigenfunction $\phi_0(x) = 1$ with the eigenvalue $0$, and the other Fourier basis functions remain invariant to the double mean centering since they are perpendicular to $\phi_0$.
Thus, as expected from Proposition~\ref{prop:S1n}, the MDS embedding $\gamma$ of $S^1$ is \begin{align*}\gamma(\theta)=  \sqrt{2}(\cos \theta, \sin \theta, \tfrac{1}{3}\cos 3\theta, \tfrac{1}{3}\sin 3\theta, \tfrac{1}{5}\cos 5\theta, \tfrac{1}{5}\sin 5\theta, \ldots ) \in \ell^2,
\end{align*} where the $\sqrt{2}$ is a normalization factor we picked up moving from a complex to a real eigendecomposition.

A couple of observations:
\begin{enumerate}
\item Applying the $\ell^2$ Euclidean distance formula to the image of $\gamma$ shows that for all $\theta \in S^1$,
\[\norm{\gamma(\theta)}_{\ell^2}^2=2\sum_{k\text{ odd}} \frac{1}{k^2} = \frac{\pi^2}{4}.\]
That is, the MDS embedding lies on an infinite-dimensional sphere of radius $\frac{\pi}{2}$ in $\ell^2$.

\item The $\ell^2$ distance between $\gamma(\theta_1)$ and $\gamma(\theta_2)$ gives an approximation of the arc-length distance between angles $\theta_1$ and $\theta_2$:
\[(\theta_1 -\theta _2)^2 \approx \norm{\gamma(\theta_1)-\gamma(\theta_2)}^2_{\ell^2} =4 \sum_{k\text{ odd}} \frac{\left(1 - \cos \left(k(\theta_1 - \theta_2) \right) \right)}{k^2}.\]
We leave it to the reader to verify that the expression above constitutes the odd modes in the Fourier series expansion of the periodic function $(\theta_1 - \theta_2)^2$.
In fact, the error of MDS comes precisely from the even modes: \begin{align*}
(\theta_1 -\theta _2)^2&=\norm{\gamma(\theta_1)-\gamma(\theta_2)}^2_{\ell^2}-\left( 4\sum_{k\text{ even}}\frac{1-\cos(k(\theta_1-\theta_2))}{k^2} \right).
\end{align*}
\end{enumerate}

For this example, the issue of convergence of MDS on finite samples to MDS on the manifold is intuitively clear: the discrete Fourier modes converge (pointwise on the sample points) to the Fourier basis $\phi(\theta)=e^{ik\theta}$.
However, in general the issue of convergence is not as straightforward.
In the next section of the paper we survey results on convergence.

The MDS embeddings of the geodesic circle are closely related to~\cite{von1941fourier}, which was written prior to the invention of MDS.
In~\cite[Theorem~1]{von1941fourier}, von Neumann and Schoenberg describe (roughly speaking) which metrics on the circle one can isometrically embed into the Hilbert space $\ell^2$.
The geodesic metric on the circle is not one of these metrics.
However, the MDS embedding of the geodesic circle into $\ell^2$ must produce a metric on $S^1$ which is of the form described in~\cite[Theorem~1]{von1941fourier}.
See also~\cite[Section~5]{wilson1935certain} and~\cite{blumenthal1970theory,bogomolny2003spectral,dattorro2010convex}.

\section{Convergence of MDS}
\label{sec: convergence}
We saw in the prior section how MDS on an evenly-spaced sample from the geodesic circle generalizes to the MDS integral operator on the entire circle.
In this section, we address convergence questions for MDS more generally.
Convergence is well-understood when each metric space has the same finite number of points~\cite{sibson1979studies}, but we are also interested in convergence when the number of points varies and is possibly infinite.

\subsection{Robustness of MDS with Respect to Perturbations}\label{ss: robustness}

\begin{figure}[htb]
   \includegraphics[width=0.3\textwidth]{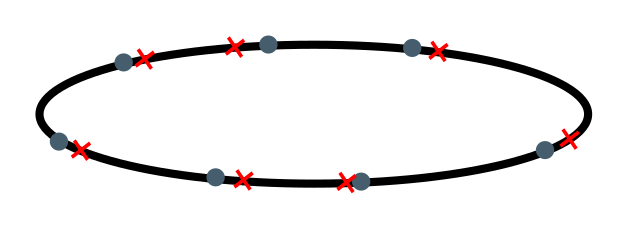}
    
    \caption{Perturbation of the given dissimilarities.}
    \label{fig:measure_Sibson}

\end{figure}

In a series of papers~\cite{sibson1978studies,sibson1979studies,sibson1981studies}, the authors consider the robustness of multidimensional scaling with respect to perturbations of the underlying dissimilarity or distance matrix, as illustrated in Figure~\ref{fig:measure_Sibson}.
In particular,~\cite{sibson1979studies} gives quantitative control over the perturbation of the eigenvalues and vectors determining an MDS embedding in terms of the perturbations of the dissimilarities.
These results build upon the fact that if $\lambda$ and $v$ are a simple (i.e., non-repeated) eigenvalue and eigenvector of an $n\times n$ matrix $\bB$, then one can control the change in $\lambda$ and $v$ upon a small symmetric perturbation of the entries in $\bB$.

Sibson's perturbation analysis shows that if one has a converging sequence of $n\times n$ dissimilarity matrices, then the corresponding MDS embeddings of $n$ points into Euclidean space also converge.
In the following sections, we consider the  convergence of MDS when the number of points is not fixed.
Indeed, we study the convergence of MDS when the number of points is finite but tending to infinity, and alternatively also when the number of points is infinite at each stage in a converging sequence of metric measure spaces.

\subsection{Convergence of MDS by the Law of Large Numbers}\label{ss: convergence random}

\begin{figure}[htb]
\centering
\begin{minipage}{0.7\linewidth}
\centering
\includegraphics[width=0.8\textwidth]{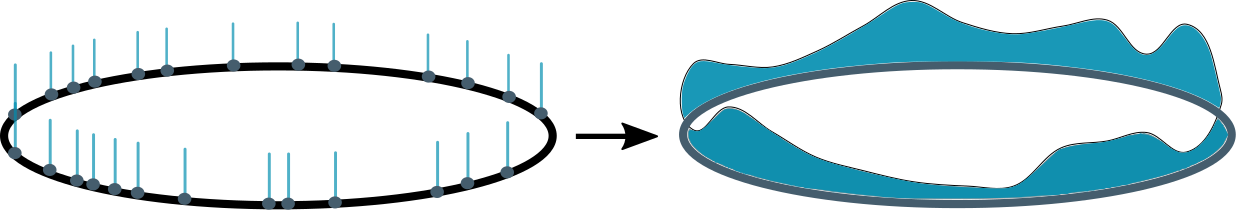}
\end{minipage}
\caption{Illustration of a notion of convergence of measures discussed in~\cite{bengio2004learning}.}
\label{fig: bengio_finiteV2}
\end{figure}

Whereas Sibson's perturbation analysis was for MDS on a fixed number of points, we now survey results
on the convergence of MDS when $n$ points $\{x_1,\ldots,x_n\}$ are sampled from a metric space according to a probability measure $\mu$, in the limit as $n\to\infty$, i.e.\ when more and more points are sampled.
In~\cite{bengio2004learning}, Bengio et al.\ study converging measures which are averages of Dirac delta functions, namely $\mu_n=\frac{1}{n}\sum_{i=1}^n\delta_{x_i}$, with all $n$ of the random points $x_i$ weighted equally (see Figure~\ref{fig: bengio_finiteV2}).
Unsurprisingly, these results rely on the law of large numbers.

Consider a data set $X_n=\{x_1, \ldots, x_n\}$ sampled independent and identically distributed (i.i.d.) from an unknown probability measure $\mu$ on $X$.
To generalize MDS, Bengio et al.\ define a corresponding data-dependent kernel that generalizes the mean centering matrix $B$ (as defined in Section~\ref{ss:infinite-mds}).
Consequently, they study the convergence of eigenvalues and eigenfuctions of the integral operator associated to the kernel as the number of sampled points increases, and they show the convergence of the MDS embeddings under desirable conditions.
They use a fundamental result on the convergence of eigenvalues of this type of integral operator from~\cite{koltchinskii2000random}.

\subsection{Convergence of MDS for Arbitrary Measures}\label{ss: convergence deterministic arbitrary}

In~\cite{kassab2019multidimensional}, we reprove the results of the previous section under a a different setting which is more general in the sense that we allow for an arbitrary sequence of convergent measures, but which is easier in the sense that this sequence is fixed (i.e.\ deterministic, not random).

\begin{figure}[htb]
\centering
\begin{minipage}{0.7\linewidth}
\centering
\includegraphics[width=0.8\textwidth]{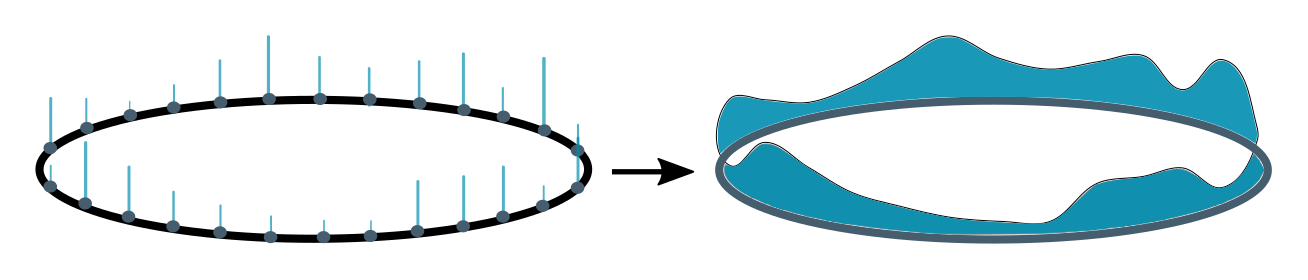}
 \subcaption{Convergence of arbitrary measures with finite support.}
\label{fig:sub:measure_deterministic_finite}
\end{minipage}

\begin{minipage}{0.7\linewidth}
\centering
\includegraphics[width=0.8\textwidth]{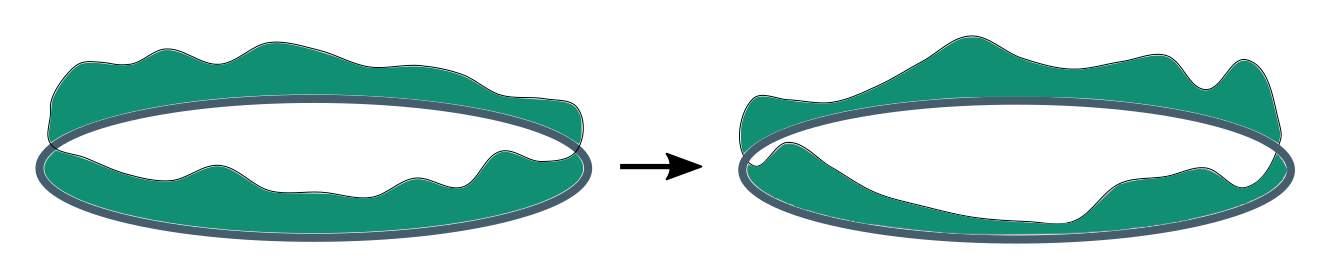}
 \subcaption{Convergence of arbitrary measures with infinite support.}
\label{fig:sub:measure_deterministic_infinite}
\end{minipage}
      \caption{Illustration of convergence (in total variation) of arbitrary measures.}
    \label{fig: arbitrary measures}
\end{figure}

Indeed, let $X$ be a compact metric space.
Suppose $\mu_n$ is an arbitrary sequence of probability measures on $X_n$ for all $n\in \N$, such that $\mu_n$ converges to $\mu$ in total variation as $n\to\infty$.
Roughly speaking, this notion of convergence of measures implies the uniform convergence of integrals against bounded measurable functions.
For example, a measure $\mu_n=\sum_{i=1}^n\lambda_i\delta_{x_i}$ in this sequence may again be a sum of Dirac delta functions, although now the weights $\lambda_i>0$ (with $\sum_i\lambda_i=1$) need not be identically equal to $\frac{1}{n}$ (Figure~\ref{fig:sub:measure_deterministic_finite}).
Much more generally, the support of any $\mu_n$ is now allowed to be infinite, as illustrated in Figure~\ref{fig:sub:measure_deterministic_infinite}.
Following~\cite{bengio2004learning,koltchinskii2000random}, we give some first results towards showing that the MDS embeddings of $(X, d, \mu_n)$ converge to the MDS embedding of $(X, d, \mu)$~\cite{kassab2019multidimensional}.
We similarily define a data-dependent kernel that generalizes the mean centering matrix $B$ (as defined in Section~\ref{ss:infinite-mds}).
It is important to note that these kernels depend on the measure on the space.
We again show convergence of eigenfunctions and consequently of MDS embeddings.

\subsection{Convergence of MDS with Respect to Gromov--Wasserstein Distance}\label{sec: convergence GW}
\begin{figure}[htb]
\centering
\begin{minipage}{0.7\linewidth}
\centering
\includegraphics[width=0.8\textwidth]{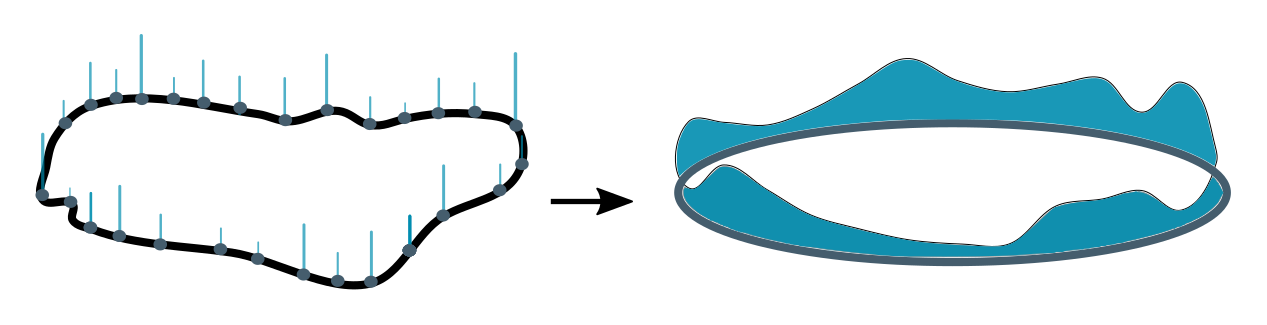}
\subcaption{Convergence of mm-spaces equipped with measures of finite support.}
\label{fig:sub:measure_GW_finite}
\end{minipage}

\begin{minipage}{0.7\linewidth}
\centering
\includegraphics[width=0.8\textwidth]{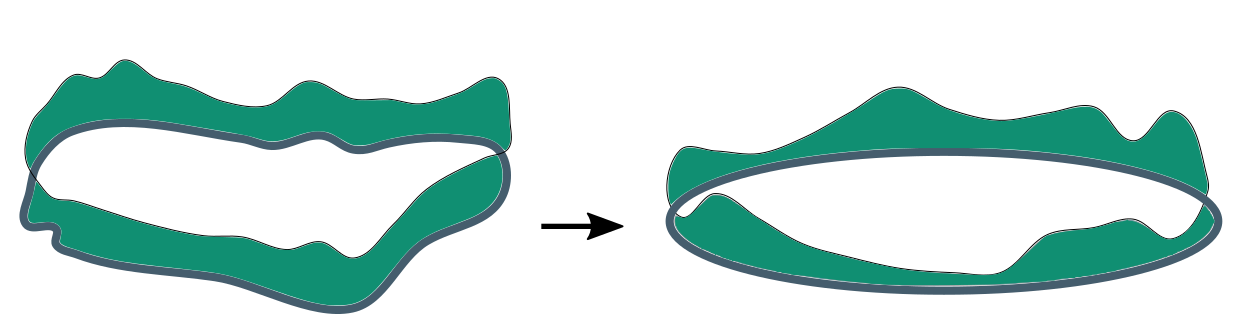}
\subcaption{Convergence of mm-spaces equipped with measures of infinite support.}
\label{fig:sub:measure_GW_infinite}
\end{minipage}

\caption{Illustration of Gromov--Wasserstein convergence of arbitrary metric measure spaces (mm-spaces).}
\label{fig: GW convergence} 
\end{figure}

We now consider the more general setting in which $(X_n, d_n, \mu_n)$ is an arbitrary sequence of metric measure spaces, converging to $(X, d, \mu)$ in the Gromov--Wasserstein distance, as illustrated in Figure~\ref{fig:sub:measure_GW_finite} for the finite case and Figure~\ref{fig:sub:measure_GW_infinite} for the infinite case.
We remark that $X_n$ need to no longer equal $X$, nor even be a subset of $X$.
Indeed, the metric $d_n$ on $X_n$ is allowed to be different from the metric $d$ on $X$.
Sections~\ref{ss: convergence random} and \ref{ss: convergence deterministic arbitrary} are the particular case when $(X_n,d_n)=(X,d)$ for all $n$, and the measures $\mu_n$ are converging to $\mu$.
We now want to consider the case where metric $d_n$ need no longer be equal to $d$.

The Wasserstein (or Kantorovich--Rubinstein) metric is a distance function defined between probability distributions on a given metric space $X$.
Intuitively, if each distribution is viewed as a unit amount of ``dirt" piled on $X$, the distance between two distributions is the minimum amount of work required to transform one pile of dirt into the other.
More generally, the Gromov--Wasserstein distance between metric measure spaces takes into account not only the variation in measures, but also the variation in metrics between these spaces.
Applications of the notion of Gromov--Wasserstein distance arise in shape and data analysis~\cite{memoli2014gromov}.

\begin{conjecture}\label{conj: GW}
Let $(X_n, d_n, \mu_n)$ for $n\in\N$ be a sequence of metric measure spaces that converges to $(X, d, \mu)$ in the Gromov--Wasserstein distance.
Then the MDS embeddings converge.
\end{conjecture}

\begin{question}\label{ques: GW}
Are there other notions of convergence of a sequence of arbitrary (possibly infinite) metric measure spaces $(X_n, d_n, \mu_n)$ to a limiting metric measure space $(X, d, \mu)$ that would imply that the MDS embeddings converge in some sense?
We remark that one might naturally try to break this into two steps: first analyze which notions of convergence $(X_n, d_n, \mu_n) \to (X, d, \mu)$ imply that the corresponding operators converge, and then analyze which notions of convergence on the operators imply that their eigendecompositions and MDS embeddings converge.
\end{question}

\section{Conclusion}
\label{chap: conclusion}
MDS is concerned with problem of mapping the objects $x_1, \ldots, x_n$ to a configuration (or embedding) of points $f(x_1), \ldots, f(x_n)$ in $\R^m$ in such a way that the given dissimilarities $d_{ij}$ are well-approximated by the Euclidean distances between $f(x_i) $ and $f(x_j)$.
We study a notion of MDS on metric measure spaces, which can be simply thought of as spaces of (possibly infinitely many) points equipped with some probability measure.
We explain how MDS generalizes to metric measure spaces.
Furthermore, we describe in a self-contained fashion an infinite analogue to the classical MDS algorithm.
Indeed, classical multidimensional scaling can be described either as a strain-minimization problem, or as a linear algebra algorithm involving eigenvalues and eigenvectors.
We describe how to generalize both of these formulations to metric measure spaces.
We show that this infinite analogue minimizes a strain function similar to the strain function of classical MDS.

As a motivating example for convergence of MDS, we consider the MDS embeddings of the circle equipped with the (non-Euclidean) geodesic metric.
By using the known eigendecomposition of circulant matrices, we identify the MDS embeddings of evenly-spaced points from the geodesic circle into $\R^m$, for all $m$.
Indeed, the MDS embeddings of the geodesic circle are closely related to~\cite{von1941fourier}, which was written prior to the invention of MDS.

Lastly, we address convergence questions for MDS.
Convergence is understood when each metric space in the sequence has the same finite number of points, or when each metric space has a finite number of points tending to infinity.
We are also interested in notions of convergence when each metric space in the sequence has an arbitrary (possibly infinite) number of points.
For instance, if a sequence of metric measure spaces converges to a fixed metric measure space $X$, then in what sense do the MDS embeddings of these spaces converge to the MDS embedding of $X$?

Several questions remain open.
In particular, we would like to have a better understanding of the convergence of MDS under the most unrestrictive assumptions of a sequence of arbitrary (possibly infinite) metric measure spaces converging to a fixed metric measure space.
Is there a version that holds under convergence in the Gromov--Wasserstein distance, which that allows for distortion of both the metric and the measure simultaneously (see Conjecture~
\ref{conj: GW} and Question~\ref{ques: GW})?
Despite all of the work that has been done on MDS by a wide variety of authors, many interesting questions remain open (at least to us).
For example, consider the MDS embeddings of the $n$-sphere for $n\ge 2$.

\begin{question}
What are the MDS embeddings of the $n$-sphere $S^n$, equipped with the geodesic metric, into Euclidean space $\R^m$?
\end{question}

To our knowledge, the MDS embeddings of $S^n$ into $\R^m$ are not understood for all positive integers $m$ except in the case of the circle, when $n=1$.
The above question is also interesting, even in the case of the circle, when the $n$-sphere is not equipped with the uniform measure.
As a specific case, what is the MDS embedding of $S^1$ into $\R^m$ when the measure is not uniform on all of $S^1$, but instead (for example) uniform with mass $\frac{2}{3}$ on the northern hemisphere, and uniform with mass $\frac{1}{3}$ on the southern hemisphere?

We note the work of Blumstein and Kvinge \cite{blumstein2018letting}, where a finite group representation theoretic perspective on MDS is employed.
Adapting these techniques to the analytical setting of compact Lie groups may prove fruitful for the case of infinite MDS on higher dimensional spheres.

We also note the work \cite{blumstein2019pseudo}, where the theory of an MDS embedding into pseudo Euclidean space is developed.
In this setting, both positive and negative eigenvalues are used to create an embedding.
In the example of embedding $S^1$, positive and negative eigenvalues occur in a one-to-one fashion.
We wonder about the significance of the full spectrum of eigenvalues for the higher dimensional spheres.

\section{Acknowledgements}
We would like to thank Bailey Fosdick, Michael Kirby, Henry Kvinge, Facundo M{\'e}moli, Louis Scharf, the students in Michael Kirby's Spring 2018 class, and the Pattern Analysis Laboratory at Colorado State University for their helpful conversations and support throughout this project.

\bibliographystyle{plain}

\bibliography{MultidimensionalScalingOnMetricMeasureSpaces.bib}

\end{document}